\documentclass{amsart}

\usepackage{amsmath,amsthm,amsfonts,amscd,amssymb}
\usepackage{url}
\usepackage{comment}

\newtheorem{thm}{Theorem}[section]
\newtheorem{cor}[thm]{Corollary}
\newtheorem{lem}[thm]{Lemma}
\newtheorem{prop}[thm]{Proposition}

\theoremstyle{definition}

\theoremstyle{remark}
\newtheorem{rem}{Remark}[section]

\begin{document}

\title{Positive semigroups in lattices and totally real number fields}
\author{Lenny Fukshansky}\thanks{Fukshansky was partially supported by the Simons Foundation grant \#519058}
\author{Siki Wang}

\address{Department of Mathematics, 850 Columbia Avenue, Claremont McKenna College, Claremont, CA 91711, USA}
\email{lenny@cmc.edu}
\address{Department of Mathematics, 850 Columbia Avenue, Claremont McKenna College, Claremont, CA 91711, USA}
\email{swang21@students.claremontmckenna.edu}


\begin{abstract}
Let $L$ be a full-rank lattice in $\mathbb R^d$ and write $L^+$ for the semigroup of all vectors with nonnegative coordinates in $L$. We call a basis $X$ for $L$ positive if it is contained in $L^+$. There are infinitely many such bases, and each of them spans a conical semigroup $S(X)$ consisting of all nonnegative integer linear combinations of the vectors of $X$. Such $S(X)$ is a sub-semigroup of $L^+$, and we investigate the distribution of the gaps of $S(X)$ in $L^+$, i.e. the points in $L^+ \setminus S(X)$. We describe some basic properties and counting estimates for these gaps. Our main focus is on the restrictive successive minima of $L^+$ and of $L^+ \setminus S(X)$, for which we produce bounds in the spirit of Minkowski's successive minima theorem and its recent generalizations. We apply these results to obtain analogous bounds for the successive minima with respect to Weil height of totally positive sub-semigroups of ideals in totally real number fields. 
\end{abstract}

\maketitle

\def\A{{\mathcal A}}
\def\AA{{\mathfrak A}}
\def\B{{\mathcal B}}
\def\C{{\mathcal C}}
\def\D{{\mathcal D}}
\def\EE{{\mathfrak E}}
\def\E{{\mathcal E}}
\def\F{{\mathcal F}}
\def\x{{\mathcal H}}
\def\I{{\mathcal I}}
\def\II{{\mathfrak I}}
\def\J{{\mathcal J}}
\def\K{{\mathcal K}}
\def\kk{{\mathfrak K}}
\def\L{{\mathcal L}}
\def\LL{{\mathfrak L}}
\def\M{{\mathcal M}}
\def\mm{{\mathfrak m}}
\def\MM{{\mathfrak M}}
\def\N{{\mathcal N}}
\def\O{{\mathcal O}}
\def\OO{{\mathfrak O}}
\def\PP{{\mathfrak P}}
\def\R{{\mathcal R}}
\def\PNR{{\mathcal P_N(\real)}}
\def\PMNR{{\mathcal P^M_N(\real)}}
\def\PdNR{{\mathcal P^d_N(\real)}}
\def\s{{\mathcal S}}
\def\V{{\mathcal V}}
\def\X{{\mathcal X}}
\def\Y{{\mathcal Y}}
\def\Z{{\mathcal Z}}
\def\H{{\mathcal H}}
\def\BB{{\mathbb B}}
\def\cee{{\mathbb C}}
\def\Nn{{\mathbb N}}
\def\pee{{\mathbb P}}
\def\que{{\mathbb Q}}
\def\QQ{{\mathbb Q}}
\def\real{{\mathbb R}}
\def\RR{{\mathbb R}}
\def\zed{{\mathbb Z}}
\def\ZZ{{\mathbb Z}}
\def\aaa{{\mathbb A}}
\def\ff{{\mathbb F}}
\def\VV{{\mathbb V}}
\def\kk{{\mathfrak K}}
\def\qbar{{\overline{\mathbb Q}}}
\def\kbar{{\overline{K}}}
\def\ybar{{\overline{Y}}}
\def\kkbar{{\overline{\mathfrak K}}}
\def\ubar{{\overline{U}}}
\def\eps{{\varepsilon}}
\def\ahat{{\hat \alpha}}
\def\bhat{{\hat \beta}}
\def\gt{{\tilde \gamma}}
\def\h{{\tfrac12}}
\def\be{{\boldsymbol e}}
\def\bei{{\boldsymbol e_i}}
\def\bc{{\boldsymbol c}}
\def\bm{{\boldsymbol m}}
\def\bk{{\boldsymbol k}}
\def\bi{{\boldsymbol i}}
\def\bl{{\boldsymbol l}}
\def\bq{{\boldsymbol q}}
\def\bu{{\boldsymbol u}}
\def\bt{{\boldsymbol t}}
\def\bs{{\boldsymbol s}}
\def\bv{{\boldsymbol v}}
\def\bw{{\boldsymbol w}}
\def\bx{{\boldsymbol x}}
\def\bX{{\boldsymbol X}}
\def\bz{{\boldsymbol z}}
\def\bwy{{\boldsymbol y}}
\def\bY{{\boldsymbol Y}}
\def\bL{{\boldsymbol L}}
\def\ba{{\boldsymbol a}}
\def\baa{{\boldsymbol \alpha}}
\def\bb{{\boldsymbol b}}
\def\bbb{{\boldsymbol \beta}}
\def\bet{{\boldsymbol\eta}}
\def\bxi{{\boldsymbol\xi}}
\def\bo{{\boldsymbol 0}}
\def\bone{{\boldsymbol 1}}
\def\bol{{\boldsymbol 1}_L}
\def\ep{\varepsilon}
\def\p{\boldsymbol\varphi}
\def\q{\boldsymbol\psi}
\def\rank{\operatorname{rank}}
\def\aut{\operatorname{Aut}}
\def\lcm{\operatorname{lcm}}
\def\sgn{\operatorname{sgn}}
\def\spn{\operatorname{span}}
\def\md{\operatorname{mod}}
\def\Norm{\operatorname{Norm}}
\def\dim{\operatorname{dim}}
\def\det{\operatorname{det}}
\def\Vol{\operatorname{Vol}}
\def\Area{\operatorname{Area}}
\def\rk{\operatorname{rk}}
\def\ord{\operatorname{ord}}
\def\ker{\operatorname{ker}}
\def\div{\operatorname{div}}
\def\Gal{\operatorname{Gal}}
\def\GL{\operatorname{GL}}
\def\SNR{\operatorname{SNR}}
\def\WR{\operatorname{WR}}
\def\scg{\operatorname{\left< \Gamma \right>}}
\def\swrh{\operatorname{Sim_{WR}(\Lambda_h)}}
\def\ch{\operatorname{C_h}}
\def\cht{\operatorname{C_h(\theta)}}
\def\scgt{\operatorname{\left< \Gamma_{\theta} \right>}}
\def\scgmn{\operatorname{\left< \Gamma_{m,n} \right>}}
\def\gat{\operatorname{\Omega_{\theta}}}
\def\GG{\operatorname{G_{\Lambda_{\ba}}}}
\def\Gg{\operatorname{G}}
\def\Sg{\operatorname{Sg}}
\def\Int{\operatorname{int}}
\def\disc{\operatorname{disc}}
\def\pr{\operatorname{pr}}
\def\Tr{\operatorname{Tr}}
\def\supp{\operatorname{supp}}
\def\d{\partial}
\def\adj{\operatorname{adj}}

\section{Introduction}
\label{intro}

Let $d \geq 2$ and $L \subset \real^d$ be a lattice of full rank, and let us write
$$\real^d_{\geq 0} = \left\{ \bx \in \real^d : x_i \geq 0\ \forall\ 1 \leq i \leq d \right\}$$
for the positive orthant of the Euclidean space $\real^d$ and $\real^d_{> 0}$ for its interior. Define
$$L^+ = L \cap \real^d_{\geq 0},$$
then $L^+$ is an additive monoid in $L$. If $X = \{\bx_1,\dots,\bx_n\}$ is a basis for $L$ contained in $L^+$, which we refer to as a {\it positive basis} for $L$, we can write 
\begin{equation}
\label{X_matrix}
\X = (\bx_1\ \dots \bx_d)
\end{equation}
for the corresponding $d \times d$ positive basis matrix, so $L = \X\zed^d$. A choice of such $X$ is always possible, since every lattice contains infinitely many positive bases (Lemma~\ref{basis+}). Define a submonoid of $L^+$
$$S(X) = \left\{ \sum_{i=1}^n a_i \bx_i : a_1,\dots,a_n \in \zed_{\geq 0} \right\} = \X \zed^d_{\geq 0},$$
as well as the positive cone spanned by $X$
$$\C(X) = \left\{ \sum_{i=1}^n a_i \bx_i : a_1,\dots,a_n \in \real_{\geq 0} \right\} = \X \real^d_{\geq 0}.$$
Define the set of {\it gaps} of $S(X)$ in $L^+$ to be $G(X) = L^+ \setminus S(X)$. We provide a basic characterization of gaps in Section~\ref{lattices}, in particular showing that $S(X) = L \cap \C(X)$, meaning that no gaps can be contained in the positive cone $\C(X)$. Hence $G(X)$ can be described as $L^+ \setminus \C(X)$, which is an infinite set unless the basis $X$ is orthogonal: in the orthogonal case, $L^+ = S(X)$. Furthermore, if $X$ is not orthogonal, then $L^+$ is not contained in any set of the form
$$\left\{ \sum_{i=1}^d a_i \bx_i : a_1,\dots,a_d \in \Z \right\},$$
where $\Z \subset \zed$ is a proper subset (Lemma~\ref{coeff_extend}). In other words, the gaps cannot all be ``filled" by partially expanding the set of coefficients to include some (but not all) negative integers. These observations motivate further investigation of the distribution of the gaps. First, we can obtain an asymptotic estimate for the number of gaps of bounded norm of the positive semigroup $S(X)$ in $L^+$ as the bound tends to infinity. For a positive real number $t$ and a set $S \subseteq \real^d$, define the counting function
$$\N(S,t) = \left| \left\{ \bwy \in S : \|\bwy\| \leq t \right\} \right|,$$
where $\|\ \|$ stands for the usual Euclidean norm on~$\real^d$.

\begin{prop} \label{count_gaps} Let $L \subset \real^d$ be a lattice of full rank and $X$ a positive basis for~$L$. Let $\nu(X)$ stand for the measure of the solid angle of the cone $\C(X)$ and write $\omega_d$ for the volume of a unit ball in~$\real^d$. As $t \to \infty$, 
\begin{equation}
\label{gaps_cnt}
\N(G(X),t) = \left( \frac{\omega_d (1 - \nu(X) 2^d)}{2^d \det L} \right) t^{d} + O(t^{d-1}).
\end{equation}
\end{prop}

The measure $\nu(X)$ of the solid angle of the cone $\C(X)$ present in the bound of~\eqref{gaps_cnt} above is defined as the proportion of the volume of the unit ball centered at the origin which is cut-out by this cone, i.e.
$$\nu(X) = \frac{\Vol_d (\C(X) \cap B_d(1))}{\omega_d},$$
where
$$B_d(t) = \left\{ \bx \in \real^d : \|\bx\| \leq t \right\}$$
is a ball of radius $t > 0$ centered at the origin in~$\real^d$. Since $\C(X) \subseteq \real^d_{\geq 0}$, the solid angle $\nu(X)$ is no bigger than $1/2^d$, meaning that the constant in the main term of~\eqref{gaps_cnt} is nonnegative, and in fact positive whenever the basis $X$ is not orthogonal. While there are known computational formulas for solid angles of cones (see, for instance, Theorem~2.2 of~\cite{ribando}) they are complicated and technical, which is why we choose not to include them here. We prove Proposition~\ref{count_gaps} at the end of Section~\ref{lattices}.
\medskip

Our main results on lattices target the so-called restricted successive minima with respect to our semigroup $S(X)$, naturally complementing previous work on points of a lattice outside of a finite union of sublattices and hypersurfaces \cite{lf_2006},~\cite{henk_thiel}. We write $|\ |$ for the sup-norm on $\real^d$, and for a real number $t > 0$, let
$$C_d(t) = \left\{ \bx \in \real^d : |\bx| \leq t \right\}$$
be a cube of radius $t$ centered at the origin in~$\real^d$. Let 
$$\mu(L) = \min \left\{ t \in \real_{> 0} : B_d(t) + L = \real^d \right\}$$
be the inhomogeneous minimum (also called the covering radius) of $L$. We will define three different sets of successive minima with respect to the cube $C_d(1)$. First, we write $0 < \lambda_1 \leq \dots \leq \lambda_d$ for the usual successive minima of $L$ with respect to $C_d(1)$, i.e.
$$\lambda_i = \min \left\{ t \in \real_{> 0} : \dim_{\real} \spn_{\real} \left( L \cap C_d(t) \right) \geq i \right\}.$$
We also define the successive minima of $L^+$ with respect to $C_d(1)$ to be
$$\lambda_i(L^+) := \min \left\{ t \in \real_{> 0} : \dim_{\real} \spn_{\real} \left( L^+ \cap C_d(t) \right) \geq i \right\}.$$
Finally, for a positive basis $X$ of $L$ we define the restricted successive minima of $L^+ \setminus S(X)$ with respect to $C_d(1)$ in the spirit of Henk and Thiel's definition in~\cite{henk_thiel}: for each $1 \leq i \leq d$, let
$$\lambda_i(L^+,X) := \min \left\{ t \in \real_{> 0} : \dim_{\real} \spn_{\real} \left( G(X) \cap C_d(t) \right) \geq i \right\}.$$
In other words, $\lambda_i(L^+)$ (respectively, $\lambda_i(L^+,X)$) is the minimal $t$ such that there exist $i$ linearly independent vectors in $L^+$ (respectively, gaps of $S(X)$ in $L^+$) with sup-norm no bigger than~$t$, so
$$0 < \lambda_1(L^+) \leq \dots \leq \lambda_d(L^+),\ 0 < \lambda_1(L^+,X) \leq \dots \leq \lambda_d(L^+,X).$$
We obtain the following bounds on these successive minima.

\begin{thm} \label{small_gap} Let $L \subset \real^d$ be a lattice of full rank. Then 
\begin{equation}
\label{smL+}
\lambda_1(L^+) \leq 2\mu(L)+1,\ \lambda_i(L^+) \leq 2\lambda_i (\mu(L) + 1)\ \forall\ 2 \leq i \leq d.
\end{equation}
On the other hand, for each non-orthogonal positive basis $X = \{ \bx_1,\dots,\bx_d \}$ of $L$ so that no $d-1$ elements of $X$ lie in a coordinate hyperplane, we have
$$\lambda_d(L^+,X) \leq \max_{1 \leq i \leq d} \max_{1 \leq m \leq d}  \left\{ \left( \max_{1 \leq k \leq d} \left[ \frac{x_{ik}}{\sum_{j=1, j \neq i}^d x_{jk}} \right] + 1 \right) \sum_{j=1, j \neq i}^d x_{jm} - x_{im} \right\},$$
where $[\ ]$ stands for the integer part function.
\end{thm}

Notice that this result can be viewed as a version of Minkowski's Successive Minima Theorem with restrictions; in particular, explicit bounds for the successive minima~$\lambda_i(L^+)$ can be deduced from~\eqref{smL+} with the help of Minkowski's theorem (see, for instance~\cite{lek}, Section~9.1, Theorem~1) and Jarnik's inequality on inhomogeneous minimum (see, for instance~\cite{lek}, Section~13.2, Theorem~1). Indeed, it is known that
\begin{equation}
\label{mink_jarnik}
\prod_{i=1}^d \lambda_i \leq \det L,\ \mu(L) \leq \frac{\sqrt{d}}{2} \sum_{i=1}^d \lambda_i.
\end{equation}
The dimensional constant $\sqrt{d}$ appears in our upper bound on $\mu(L)$ here because we are comparing the inhomogeneous minimum with respect to the Euclidean ball $B_d(1)$ to the successive minima with respect to the cube $C_d(1)$. We prove Theorem~\ref{small_gap} in Section~\ref{sm_section}. 

In fact, one can generalize the discussion above by considering, instead of $L^+$, the intersection of $L$ with a general convex cone spanned by a basis in $\real^d$, not just with $\real^d_{\geq 0}$. Such generalizations are not difficult to obtain from our results. Say, $Y = \{ \bwy_1,\dots,\bwy_d \}$ is a collection of linear independent vectors in $\real^d$, and $L(Y) = L \cap \C(Y)$. If we write $\Y$ for the $d \times d$ matrix $(\bwy_1\ \dots\ \bwy_d)$ and let $M = \Y^{-1} L$, then $L(Y) = \Y M^+$. Thus analogues of the above observations, including Proposition~\ref{count_gaps} and Theorem~\ref{small_gap} can we obtained for $M^+$ and then ``transferred" to $L(Y)$ by applying the linear transformation~$\Y$; we present these generalizations in Corollary~\ref{gen_count} and Corollary~\ref{gen_small}, respectively. We chose to focus specifically on the positive cone $\real^d_{\geq 0}$, and hence on the positive semigroups $L^+$ because of the natural connection to totally real algebraic number fields, which we discuss next.

Let $K$ be a totally real number field of degree $d$ over $\que$, and $\sigma_1,\dots,\sigma_d : K \to \real$ be the embeddings of $K$. We write $\Nn_K$ for the field norm on $K$ and $\Delta_K$ for the discriminant of $K$. Let $\O_K$ be the ring of integers of $K$, then $\O_K$, as well as any ideal in $\O_K$, can be viewed as a Euclidean lattice of rank~$d$ with respect to the symmetric bilinear form
\begin{equation}
\label{trace}
\left< \alpha, \beta \right> = \Tr_K(\alpha \beta) = \sum_{i=1}^d \sigma_i(\alpha) \sigma_i(\beta),
\end{equation}
where $\Tr_K$ is the usual trace on $K$. For an ideal $I \subset \O_K$, let $I^+$ be the additive semigroup of totally positive elements in $I$, i.e.
$$I^+ = \left\{ \alpha \in I : \sigma_i(\alpha) \geq 0\ \forall\ 1 \leq i \leq d \right\},$$
and let us also write $\zed^+ = \zed \cap \O_K^+$. Then $I$ has a $\zed$-basis contained in $I^+$. Let $\bbb = \{ \beta_1,\dots,\beta_d \} \subset I^+$ be such a $\zed$-basis for $I$, which we call a positive basis. Let
$$S(\bbb) = \left\{ \sum_{i=1}^d c_i \beta_i : c_1,\dots,c_d \in \zed^+ \right\} \subseteq I^+$$ 
be the corresponding sub-semigroup, and define the set of gaps of $S(\bbb)$ in $I^+$ to be $G(\bbb) = I^+ \setminus S(\bbb)$. The basis $\bbb$ cannot be orthogonal, as we explain in Section~\ref{ideals}, where we discuss this setup in further details. Hence the set $G(\bbb)$ is infinite. With this notation, we can state a version of Theorem~\ref{small_gap} for ideals, where we use Weil height $h$ on $K$ to measure size of points; we recall the definition and standard properties of $h$ in Section~\ref{ideals}.

\begin{thm} \label{ideal_gaps} Let $I \subseteq \O_K$ be an ideal. Then there exist $\que$-linearly independent elements $s_1,\dots,s_d \in I$ such that
$$\prod_{i=1}^d h(s_i) \leq \Nn_K(I) \sqrt{ |\Delta_K| }.$$
Further, there exist $\que$-linearly independent elements $\alpha_1,\dots,\alpha_d \in I^+$ such that
$$\prod_{i=1}^d h(\alpha_i)  \leq \left( 3d \sqrt{d} \right)^d \left( \Nn_K(I) \sqrt{ |\Delta_K| } \right)^{d+1}.$$
Additionally, let $\bbb = \{ \beta_1,\dots,\beta_d \} \subset I^+$ be a positive basis for $I$ and $G(\bbb)$ the corresponding set of gaps. For each $1 \leq i \leq d$, let $\beta'_i = \sum_{j=1, j \neq i}^d \beta_j$. Then there exist $\que$-linearly independent gaps $\alpha_1,\dots,\alpha_d \in G(\bbb)$ such that
$$h(\alpha_i) \leq \left( h(\beta_i/\beta'_i)^d + 1\right) h(\beta'_i)^d,$$
for each $1 \leq i \leq d$.
\end{thm}

We prove this theorem in Section~\ref{ideals} by embedding our ideal into Euclidean space, applying Theorem~\ref{small_gap}, and then ``lifting" the results back up to the number field. We are now ready to proceed.
\bigskip

\section{Positive semigroups in lattices}
\label{lattices}

In this section, we develop the necessary basic observations for positive semigroups in general Euclidean lattices in a series of lemmas.

\begin{lem} \label{basis+} There exists a basis $\bx_1,\dots,\bx_d$ for $L$ contained in $\real^d_{> 0}$. Therefore there exist infinitely many such bases.
\end{lem}

\proof
Let $\bwy_1,\dots,\bwy_d \in L$ be a basis for $L$, and write $Y = (\bwy_1 \dots \bwy_d)$ for the corresponding $d \times d$ basis matrix, so $L = Y\zed^d$. Let  $\bx_1$ be any point in $L \cap \real^d_{> 0}$. Then
$$\bx_1 = \sum_{i=1}^d a_i \bwy_i$$
for some integer coefficients $a_1,\dots,a_d$. In fact, we can assume that $\gcd(a_1,\dots,a_d) = 1$: if not, divide these coefficients by their gcd and replace $\bx_1$ with the resulting vector, which is still in $L \cap \real^d_{> 0}$. Let $\ba = (a_1 \dots a_d)$ be the row vector of these relatively prime coefficients, then $\ba$ is extendable to a matrix $A \in \GL_n(\zed)$ (Lemma 2, p. 15, \cite{cass:geom}). Then $YA$ is another basis matrix for $L$, the first column vector of which is $\bx_1$, hence $\bx_1$ is extendable to a basis $\bx_1,\bx_2,\dots,\bx_d$ for $L$. Since all coordinates of $\bx_1$ are positive, we can now replace each $\bx_i$, $2 \leq i \leq d$ with $\bx_i + M_i \bx_1$ for an appropriately large integer $M_i$ to ensure that $\bx_i + M_i \bx_1$ also has all positive coordinates. The resulting collection
$$\bx_1, \bx_2 + M_2 \bx_1, \dots, \bx_d + M_d \bx_1$$
is again a basis for $L$, contained in $\real^d_{> 0}$. Once we have one such basis, we can obtain infinitely many simply by adding positive integer multiples of one of the vectors to the others. Further, there are infinitely many choices for the vector $\bx_1$, which is the starting point of our construction.
\endproof

\begin{rem} Since every matrix in $\GL_d(\real)$ is a basis matrix for some lattice, and change of basis is performed by right multiplication by a matrix from $\GL_d(\zed)$, we can rephrase Lemma~\ref{basis+} as follows: every orbit of $\GL_d(\real)$ under the action of $\GL_d(\zed)$ by right multiplication contains a matrix with all positive entries.
\end{rem}

\begin{lem} \label{C_vs_S} Let $X = \{ \bx_1,\dots,\bx_d \}$ be a positive basis for $L$, then $L^+ \cap \C(X) = S(X)$.
\end{lem}

\proof
Clearly $S(X) \subseteq L^+ \cap \C(X)$, so we only need to prove that $L^+ \cap \C(X) \subseteq S(X)$. Suppose not, then there exists some $\bwy \in L^+ \cap \C(X)$ which is not in $S(X)$. Since $\bwy \in L^+ \setminus S(X)$, we must have
$$\bwy = \sum_{i=1}^d a_i \bx_i,$$
where the coefficients $a_1,\dots,a_d \in \zed$ are not all nonnegative. On the other hand, since $\bwy \in \C(X)$, we must have
$$\bwy = \sum_{i=1}^d c_i \bx_i$$
for some nonnegative real coefficients $c_1,\dots,c_d$. Thus $\sum_{i=1}^d a_i \bx_i = \sum_{i=1}^d c_i \bx_i$, meaning that
$$\sum_{i=1}^d (a_i-c_i) \bx_i = \bo,$$
where not all of the coefficients $a_i - c_i$ can be equal to $0$. This contradicts the fact that the vectors $\bx_1,\dots,\bx_d$ are linearly independent in $\real^d$. Hence we have $L^+ \cap \C(X) = S(X)$.
\endproof

Therefore Lemmas~\ref{C_vs_S} implies that all the gaps of $S(X)$ in $L^+$ are outside of the cone $\C(X)$. On the other hand, there are infinitely many such gaps. Define the set of {\it primitive gaps} to be
$$G_{\pr}(X) = \left\{ \bwy \in G(X) : \bwy \neq m \bwy' \text{ for any } \bwy' \in L \text{ and integer } m > 1 \right\},$$
i.e. $G_{\pr}(X)$ is the set of gaps that are primitive points of the lattice $L$. 

\begin{lem} \label{Gpr_inf} The set $G_{\pr}(X)$ is infinite, unless the positive basis $X$ is orthogonal.
\end{lem}

\proof
Suppose $X$ is not an orthogonal basis, then $\C(X)$ is an acute cone, and so the set $\real^d_{\geq 0} \setminus \C(X)$ is unbounded. Thus its interior contains Euclidean balls of arbitrarily large radius, hence the union of all such balls must contain infinitely many primitive points of the lattice $L$. 
\endproof

In fact, all gaps are multiples of primitive gaps. 

\begin{lem} \label{gaps_L} For any $\bwy \in G(X)$ and $m \in \zed^+$, $m\bwy \in G(X)$. Further, for any $\bwy \in G(X)$ there exists $\bz \in G_{\pr}(X)$ and $m \in \zed^+$ such that $\bwy = m\bz$, i.e.
$$G(X) = \left\{ m \bz : \bz \in G_{\pr}(X),\ m \in \zed^+ \right\}.$$
\end{lem}

\proof
Suppose that $\bwy$ and $\bz$ are two vectors in $L^+$ such that $\bwy = m\bz$ for some $m \in \zed^+$. We just need to show that either both $\bwy,\bz$ are in $S(X)$, or they are both not. Indeed, suppose one of them, say, $\bwy$ is in $S(X)$ and $\bz$ is not. Then
$$\bwy = \sum_{i=1}^d a_i \bx_i,\ \bz = \sum_{i=1}^d c_i \bx_i,$$
where the coefficients $a_1,\dots,a_d \in \zed^+$ and some of the coefficients $c_1,\dots,c_d \in \zed$ are negative. Therefore
$$m \sum_{i=1}^d c_i \bx_i = \sum_{i=1}^d a_i \bx_i,$$
i.e. $\sum_{i=1}^d (mc_i - a_i) \bx_i = \bo$ with not all of the coefficients $mc_i-a_i = 0$. This contradicts linear independence of $\bx_1,\dots,\bx_d$, and thus completes the proof.
\endproof

\begin{rem} Notice that the above lemma is in fact equivalent to the observation that a ray emanating from the origin either does not intersect the cone $\C(X)$ or is contained in it.
\end{rem}

\begin{lem} \label{coeff_extend} Suppose that $X = \{ \bx_1,\dots,\bx_d \}$ is a positive non-orthogonal basis for $L$. Then $L^+$ is not contained in any set of the form
$$S(X,\Z) := \left\{ \sum_{i=1}^d a_i \bx_i : a_1,\dots,a_d \in \Z \right\},$$
where $\Z \subset \zed$ is a proper subset.
\end{lem}

\proof
Suppose that there exists some proper subset $\Z \subset \zed$ such that $L^+ \subseteq S(X,\Z)$, and let $b \in \zed \setminus \Z$. If $b \geq 0$, then by uniqueness of representation of vectors in $L$ with respect to $X$, the vector $b \bx_1 \in L^+ \setminus S(X,\Z)$. Hence we can assume that $b < 0$.

For a vector $\bx \in \real^d$, let us define its support as
$$\supp(\bx) := \{ 1 \leq i \leq d : x_i \neq 0 \}.$$
Since $X$ is a basis for $\real^d$,
$$\bigcup_{i=1}^d \supp(\bx_i) = \{1,\dots,d\}.$$
Assume that for every $1 \leq j \leq d$,
$$\supp(\bx_j) \not\subseteq \bigcup_{i=1, i \neq j}^d \supp(\bx_i),$$
then for each $1 \leq j \leq d$ there exists some $1 \leq k \leq d$ such that $x_{jk} \neq 0$ while $x_{ik} = 0$ for every $i \neq j$. Since there are $d$ vectors and $d$ coordinates, this is only possible if every vector has a unique nonzero coordinate, which implies that the vectors $\bx_1,\dots,\bx_d$ are simply multiples of the standard basis vectors. This, however, contradicts the assumption that $X$ is not an orthogonal basis. 

Therefore there exists some vector in $X$, say $\bx_1$, such that
$$\supp(\bx_1) \subseteq \bigcup_{i=2}^d \supp(\bx_i).$$
Since all the nonzero coordinates of $b\bx_1$ are negative and all the nonzero coordinates of $\bx_i$'s are positive, we can always pick positive integers $a_2,\dots,a_d$ large enough so that
$$\bz := b\bx_1 + \sum_{i=2}^d a_i\bx_i \in L^+,$$
however, again by uniqueness of representation of vectors with respect to the basis $X$, $\bz \not\in S(X,\Z)$, since $b \not\in \Z$. This finishes the proof.
\endproof
\medskip

Finally, we prove the counting estimate of Proposition~\ref{count_gaps}. Let $X$ be a positive basis for $L$ and $\X$ the corresponding basis matrix as in~\eqref{X_matrix}. We obtain counting estimates for the number of points of bounded norm in~$L^+$ and in the semigroup $S(X)$. Since $G(X) = L^+ \setminus S(X)$, we can then easily compute
\begin{equation}
\label{GX_diff}
\left| \left\{ \bwy \in G(X) : \|\bwy\| \leq t \right\} \right| = \left| \left\{ \bwy \in L^+ : \|\bwy\| \leq t \right\} \right| - \left| \left\{ \bwy \in S(X) : \|\bwy\| \leq t \right\} \right|.
\end{equation}

\begin{lem} \label{L+count} As $t \to \infty$,
$$\N(L^+,t) = \left( \frac{\omega_d}{2^d \det L} \right) t^{d} + O(t^{d-1}).$$
\end{lem}

\proof
It is a well-known fact (see, for instance, Theorem 2 on p.128 of~\cite{lang}) that for a lattice $L$ and a compact convex body $E$ in $\real^d$,
\begin{equation}
\label{total}
|tE \cap L| = \left( \frac{\Vol_d(E)}{\det L} \right) t^d + O(t^{d-1}),
\end{equation}
as $t \to \infty$. Let $E = B_d(1) \cap \real_{\geq 0}$, then $tE = B_d(t) \cap \real_{\geq 0}$ and $\Vol_d(E) = \omega_d/2^d$. Since $\N(L^+,t) = \left| \left( B_d(t) \cap \real^d_{\geq 0} \right) \cap L \right|$, the result follows from~\eqref{total}.
\endproof

\begin{lem} \label{SX_count} Let $\nu(X)$ be the measure of the solid angle of the cone $\C(X)$. Then
$$\N(S(X), t) = \left( \frac{\omega_d \nu(X)} {\det L} \right) t^d + O(t^{d-1}),$$
as $t \to \infty$.
\end{lem}

\proof
Let $E = B_d(1) \cap \C(X)$, then $tE = B_d(t) \cap \C(X)$ and $\Vol_d(E) = \omega_d \nu(X)$. Since $\N(S(X),t) = \left| \left( B_d(t) \cap \C(X) \right) \cap L \right|$, the result follows from~\eqref{total}.
\endproof

Proposition~\ref{count_gaps} now follows by combining~\eqref{GX_diff} with Lemmas~\ref{L+count} and \ref{SX_count}. We can also extend it to more general cones. Let $Y = \{ \bwy_1,\dots,\bwy_d \}$ be a collection of linearly independent vectors in~$\real^d$ and write $\Y = (\bwy_1\ \dots\ \bwy_d) \in \GL_d(\real)$ for the corresponding  matrix. Let $\C(Y) = \Y\real^d_{\geq 0}$ be the positive cone spanned by $Y$ and let $L(Y) = L \cap \C(Y)$. 

\begin{cor} \label{gen_count}  There exist infinitely many bases for $L$ contained in $L(Y)$. If $X$ is such a basis, then $\C(X) \subseteq \C(Y)$ and $S(X) = L \cap \C(X)$. The set of gaps $G_Y(X) := L(Y) \setminus S(X)$ is infinite unless $\C(X) = \C(Y)$ and
$$\N(G_Y(X),t) = \left( \frac{\omega_d (\nu(Y) - \nu(X))}{\det L} \right) t^{d} + O(t^{d-1})$$
as $t \to \infty$, where $\nu(Y)$ and $\nu(X)$ are the measures of the solid angles of the cones $\C(Y)$ and $\C(X)$, respectively.
\end{cor}

\proof
Let $M = \Y^{-1} L$, then $M \subset \real^d$ is also a full-rank lattice and $\det M = \frac{\det L}{|\det \Y|}$. Further, $\Y^{-1} \C(Y) = \real^d_{\geq 0}$ and $\Y^{-1} L(Y) = M^+$. Hence existence of infinitely many bases for $L$ contained in $L(Y)$ follows from the existence of infinitely many positive bases for $M$. If $X \subset L(Y)$ is such a basis for $L$, then $\Y^{-1}X$ is a positive basis for $M$, $S(X) = \Y S(\Y^{-1} X)$ and the set of gaps $G_Y(X)  = \Y G(\Y^{-1} X)$ is infinite, unless $\Y^{-1} X$ is orthogonal, which is equivalent to the assertion that $\C(X) = \C(Y)$.

To obtain the counting estimate on the number of gaps of bounded Euclidean norm, observe that
$$\N(G_Y(X,t)) = \N(L(Y),t) - \N(S(X),t),$$
where $\N(S(X),t)$ is given by Lemma~\ref{SX_count}. To obtain an estimate $\N(L(Y),t)$, we use~\eqref{total} with the choice of $E = B_d(1) \cap \C(Y)$, noting that $\Vol_d(E) = \omega_d \nu(Y)$.
\endproof

\bigskip

\section{Restricted successive minima}
\label{sm_section}

In this section we prove Theorem~\ref{small_gap}. Throughout this section we use notation as in Sections~\ref{intro} and~\ref{lattices} above, and let $X$ be a positive basis for a full-rank lattice $L \subset \real^d$.

\begin{lem} \label{min_L+} There exists a point $\bwy \in L$ such that
$$1 \leq y_i \leq 2\mu(L)+1,$$
for each $1 \leq i \leq d$. Hence
$$\lambda_1(L^+) \leq 2\mu(L)+1.$$
\end{lem}

\proof
For $r = \mu(L)$, let
$$\bz_r = (r+1) (1,\dots,1)^{\top} \in \real^d,$$
and consider the cube 
$$C_d(r) + \bz_r = \left\{ \bx \in \real^d : 1 \leq x_i \leq 2r+1 \right\}$$
in the positive orthant of $\real^d$. Then $B_d(r) + \bz_r$ is the ball of radius $r$ inscribed in $C_d(r) + \bz_r$. Since
$$\real^d = \bigcup_{\bx \in L} (B_d(r) + \bx),$$
the ball $B_d(r) + \bz_r$ must be covered by some translates of $B_d(r)$ by points of the lattice $L$. On the other hand, the ball $B_d(r) + \bz_r$ of radius $r$ cannot be covered by other balls of the same radius unless at least one of them has its center in $B_d(r) + \bz_r$. This means that there must exist 
$$\bwy \in L \cap (B_d(r) + \bz_r) \subset C_d(r) + \bz_r,$$
and so $1 \leq y_i \leq 2r+1 = 2\mu(L)+1$.
\endproof

\begin{lem} \label{sm_L+} For each $2 \leq i \leq d$,
$$\lambda_i(L^+) \leq 2\lambda_i (\mu(L) + 1).$$
\end{lem}

\proof
Let $\bwy$ be the vector constructed in Lemma~\ref{min_L+} above, and let $\bx_1,\dots,\bx_d$ be vectors corresponding to successive minima $\lambda_1,\dots,\lambda_d$. Then for at least one $1 \leq j \leq d$,
$$\bwy \notin \spn_{\real} \{ \bx_i : i \neq j \},$$
so for this $j$ let $I_j = \{ 1 \leq i \leq d : i \neq j \}$. Then the collection of $d$ vectors
$$\{ \bwy \} \cup \{ \lambda_i \bwy + \bx_i : i \in I_j \}$$
is linearly independent, and for each $i \in I_j$,
$$|\lambda_i \bwy + \bx_i| \leq \lambda_i |\bwy| + |\bx_i| \leq 2\lambda_i (\mu(L) + 1).$$
Further, since all coordinates of $\bwy$ are $\geq 1$, for each $1 \leq k \leq d$, the $k$-th coordinate of each such vector $\bwy + \bx_i$ is greater or equal than
$$\lambda_i + x_{ik} \geq 0,$$
so all of these vectors are in $L^+$. Finally notice that $\lambda_i(L^+) \leq |\bz|$ for each $\bz \in L^+$ such that
$$\dim_{\real}  \left\{ \lambda_1 \bwy + \bx_1, \dots, \lambda_{i-1} \bwy + \bx_{i-1}, \bz \right\} = i$$
if $i \leq j$ and
$$\dim_{\real} \left\{ \lambda_1 \bwy + \bx_1, \dots, \lambda_{j-1} \bwy + \bx_{j-1}, \bwy, \lambda_{j+1} \bwy + \bx_{j+1}, \dots, \lambda_{i-1} \bwy + \bx_{i-1}, \bz \right\} = i$$
if $i > j$. We can then take $\bz = \lambda_i \bwy + \bx_i$ for each $i \neq j$ and $\bz=\bwy$ for $i=j$. The lemma follows.
\endproof

\begin{lem} \label{gap_ht_1} Let $X = \{ \bx_1,\dots,\bx_d\}$ be a positive basis for $L$ so that no $d-1$ elements of $X$ lie in a coordinate hyperplane. There exist linearly independent vectors $\bz_1,\dots,\bz_d \in G_{\pr}(X)$ with
\begin{equation}
\label{gap_ht_1-eq}
|\bz_i| = \max_{1 \leq m \leq d}  \left\{ \left( \max_{1 \leq k \leq d} \left[ \frac{x_{ik}}{\sum_{j=1, j \neq i}^d x_{jk}} \right] + 1 \right) \sum_{j=1, j \neq i}^d x_{jm} - x_{im} \right\}.
\end{equation}
\end{lem}

\proof
It is enough to prove that there exist such vectors $\bz_i \in G(X)$ satisfying~\eqref{gap_ht_1-eq}: if, say, some $\bz_i$ is not in $G_{\pr}(X)$, then $\bz_i = m \bz'_i$ for some $m \in \zed^+$ and $\bz'_i \in G_{\pr}(X)$, so $|\bz'_i| \leq |\bz_i|$. Suppose that a vector $\bz = \sum_{i=1}^d a_i \bx_i \in L^+$ with at least one of the integer coefficients $a_i$ equal to $-1$. Since the representation of $\bz$ in terms of the basis $X$ is unique, and at least one of the coefficients in this representation is negative, $\bz$ is not in $S(X)$, thus $\bz \in G(X)$. 

First, we want to select a point $\bwy \in S(X)$ that would have all positive coordinates. Let
$$\bwy = \bx_1 + \dots + \bx_d.$$
Since all the coordinates of $\bx_i$'s are nonnegative and these vectors form a basis for $\real^d$, the sum of their  $k$-th coordinates has to be positive for each $1 \leq k \leq d$. Now, for each $1 \leq i \leq d$ define
$$\bz_i = a_i \bwy - (a_i+1) \bx_i = a_i \sum_{j=1, j \neq i}^d \bx_j - \bx_i,$$
where $a_i$ is a positive integer to be specified. In order for such $\bz_i$ to be in $G(X)$, we only need it to be in $L^+$, meaning that for each $1 \leq k \leq d$, we must have
$$a_i \sum_{j=1, j \neq i}^d x_{jk} > x_{ik},$$
so take 
$$a_i = \max_{1 \leq k \leq d} \left[ \frac{x_{ik}}{\sum_{j=1, j \neq i}^d x_{jk}} \right] + 1.$$
Notice that $\sum_{j=1, j \neq i}^d x_{jk} \neq 0$ since no $d-1$ elements of $X$ lie in a coordinate hyperplane. With this choice of $a_i$,
$$|\bz_i| = \max_{1 \leq m \leq d}  \left\{ \left( \max_{1 \leq k \leq d} \left[ \frac{x_{ik}}{\sum_{j=1, j \neq i}^d x_{jk}} \right] + 1 \right) \sum_{j=1, j \neq i}^d x_{jm} - x_{im} \right\},$$
which yields~\eqref{gap_ht_1-eq}. Finally, notice that the vectors $\bz_1,\dots,\bz_d$ are linearly independent since the vectors $\bx_1,\dots,\bx_d$ are.
\endproof

Now Theorem~\ref{small_gap} follows from Lemmas~\ref{sm_L+} and~\ref{gap_ht_1}. We can also extend it to more general cones $\C(Y) = \Y \real^d_{\geq 0}$ as in Corollary~\ref{gen_count}, where $Y = \{ \bwy_1,\dots,\bwy_d \} \subset \real^d$ is a collection of linearly independent vectors and $\Y = (\bwy_1\ \dots\ \bwy_d) \in \GL_d(\real)$ is the corresponding  matrix. Write $|\Y|$ and $|\Y^{-1}|$ for the sup-norm of matrices $\Y$ and $\Y^{-1}$, respectively, viewed as vectors in $\real^{d^2}$, and for each $1 \leq i \leq d$ define
$$\lambda_i(L(Y)) := \min \left\{ t \in \real_{> 0} : \dim_{\real} \spn_{\real} \left( L(Y) \cap C_d(t) \right) \geq i \right\},$$
and
$$\lambda_i(L(Y),X) := \min \left\{ t \in \real_{> 0} : \dim_{\real} \spn_{\real} \left( G_Y(X) \cap C_d(t) \right) \geq i \right\}.$$
These successive minima are easy to relate to the corresponding successive minima of the lattice $M = \Y^{-1} L$.

\begin{cor} \label{gen_small} Let $M = \Y^{-1} L$, then for each $1 \leq i \leq d$,
$$\frac{\lambda_i(M^+)}{d |\Y^{-1}|} \leq \lambda_i(L(Y)) \leq d |\Y| \lambda_i(M^+).$$
Additionally,
$$\lambda_d(L(Y), X) \leq d |\Y| \lambda_d(M^+, \Y^{-1} X).$$

\end{cor}

\proof
Let $M = \Y^{-1} L$, then $M^+ = L(Y)$ and a vector $\bz_i \in M^+$ corresponds to the successive minimum $\lambda_i(M^+)$ if and only if $\Y \bz_i$ corresponds to the successive minimum $\lambda_i(L(Y))$. Then
$$\lambda_i(L(Y)) = |\Y\bz_i| \leq d |\Y| |\bz_i| = d |\Y| \lambda_i(M^+),$$
as well as
$$\lambda_i(M^+) = |\Y^{-1} (\Y\bz_i)| \leq d |\Y^{-1}| |\Y\bz_i| = d |\Y^{-1}| \lambda_i(L(Y)).$$
Further, $\bz \in G_Y(X)$ corresponds to the successive minimum $\lambda_i(L(Y), X)$ if and only $\Y^{-1} \bz$ is a gap of $S(\Y^{-1} X)$ in $M^+$ corresponding to the successive minimum $\lambda_i(M^+, \Y^{-1} X)$. This observation combined with the inequalities above completes the proof.
\endproof

\bigskip

\section{Positive semigroups in number fields}
\label{ideals}

In this section, we let $K$ be a totally real number field and use the notation of Section~\ref{intro}. Investigation of totally positive semigroups in totally real number fields from the standpoint of the generalized Frobenius problem, a somewhat different perspective, has been initiated in~\cite{lf_shi}. We can view an ideal $I \subseteq \O_K$ as a lattice embedded into the Euclidean space $\real^d$ via the Minkowski embedding
$$\Sigma = (\sigma_1,\dots,\sigma_d) : K \to \real^d.$$
Indeed, the image of our ideal $I$ under $\Sigma$, $\Sigma(I)$ is a lattice of full rank in $\real^d$, call it $L_I$, and $L_I^+ = \Sigma(I^+)$. The Euclidean norm $\|\ \|$ on $L_I$ then precisely corresponds to the trace form on $K$, defined in~\eqref{trace}, i.e.
$$\|\Sigma(\alpha)\|^2 = \Tr_K(\alpha^2),$$
for any $\alpha \in I$. Further,
\begin{equation}
\label{det_norm}
\det L_I = \Nn_K(I) \sqrt{ |\Delta_K| },
\end{equation}
where $\Nn_K$ is the norm on $K$ and $\Delta_K$ the discriminant of $K$, as stated in Section~\ref{intro} (see, for instance, Lemma 2 on p.~115 of~\cite{lang}). Notice that a basis $\bbb = \{ \beta_1,\dots,\beta_d \}$ for $I$ is contained in $I^+$ if and only if its image 
$$X = \{ \bx_1,\dots,\bx_d \} := \left\{ \Sigma(\beta_1),\dots,\Sigma(\beta_d) \right\}$$
is a basis for $L_I$ contained in $L_I^+$, and $\Sigma(S(\bbb)) = S(X)$. Then Lemma~\ref{basis+} implies that $I$ has infinitely many bases contained in $I^+$. Further, for a positive basis $\bbb$ of~$I$, $\alpha$ is a gap of $S(\bbb)$ if and only if $\Sigma(\alpha)$ is a gap of $S(\Sigma(\bbb))$. 

The basis $\bbb$ is not orthogonal. Indeed, since $X \subset L_I^+$, in order for $\bbb$ to be orthogonal the vectors of $X$ must be along the coordinate axes in $\real^d$, meaning that these vectors have zero coordinates. This is not possible, since coordinates of any nonzero vector $\bwy = \Sigma(\alpha) \in L_I$ are conjugates of a nonzero element $\alpha \in I$, hence cannot be zero. Therefore the set of gaps $G(\bbb) := I^+ \setminus S(\bbb)$ is infinite.

To measure size of elements in our number field, we use Weil height. Let $M(K)$ be the set of all places of $K$, and for every $\alpha \in K$ let
$$h(\alpha) = \prod_{v \in M(K)} \max \{1,|\alpha|_v\}^{d_v/d},$$
where $d_v = [K_v:\que_v]$ is the local degree of $K$ at the place $v \in M(K)$. Notice that for each $v \mid \infty$, $d_v = 1$ since $K$ is totally real. 

\begin{lem} \label{ht_ineq} For every nonzero $\alpha \in \O_K$,
$$1 \leq h(\alpha) \leq |\Sigma(\alpha)|,$$
and for every nonzero $\alpha \in K$,
$$|\Sigma(\alpha)| \leq h(\alpha)^d.$$
\end{lem}

\proof
Let $\alpha \in \O_K$ be nonzero. For each $v \nmid \infty$, $|\alpha|_v \leq 1$, and so
$$h(\alpha) = \prod_{v \mid \infty} \max \{1, |\alpha|_v\}^{1/d} = \left( \prod_{i=1}^d \max \{1, |\sigma_i(\alpha)|\} \right)^{1/d}.$$
By the Artin-Whaples product formula~\cite{artin_whaples} combined with the arithmetic-geometric mean inequality,
$$1 = \prod_{v \in M(K)} |\alpha|^{d_v/d} \leq \prod_{i=1}^d |\sigma_i(\alpha)|^{1/d} \leq \frac{1}{d} \sum_{i=1}^d |\sigma_i(\alpha)| \leq \max_{1 \leq i \leq d} |\sigma_i(\alpha)| = |\Sigma(\alpha)|.$$
This implies that
$$h(\alpha) \leq \left( \prod_{i=1}^d \max_{1 \leq j \leq d} \{|\sigma_j(\alpha)|\} \right)^{1/d} \leq |\Sigma(\alpha)|.$$
On the other hand, for any nonzero $\alpha \in K$,
$$|\Sigma(\alpha)| = \max_{1 \leq i \leq d} |\sigma_i(\alpha)| \leq \prod_{i=1}^d \max \{1, |\sigma_i(\alpha)| \} \leq \prod_{v \in M(K)} \max \{1, |\alpha|_v^{d_v} \} = h(\alpha)^d.$$
\endproof

We can now state the bounds on the successive minima of the ideal $I$ with respect to the Weil height.

\begin{lem} \label{I_sm} There exist $\que$-linearly independent elements $s_1,\dots,s_d \in I$ such that
\begin{equation}
\label{I_sm-1}
\prod_{i=1}^d h(s_i) \leq \prod_{i=1}^d |\Sigma(s_i)| \leq \Nn_K(I) \sqrt{ |\Delta_K| }.
\end{equation}
Further, the inhomogeneous minimum of the lattice $L_I$ satisfies the inequality
\begin{equation}
\label{I_sm-2}
\mu(L_I) \leq \frac{d^{3/2}}{2} \Nn_K(I) \sqrt{ |\Delta_K| }.
\end{equation}
\end{lem}

\proof
Notice that the successive minima $\lambda_i$ of the lattice $L_I$, as we defined them in Section~\ref{intro} are the smallest sup-norms of $\que$-linearly independent vectors in $L_I$. Taking $s_1,\dots,s_d \in I$ to be a collection of elements so that $\Sigma(s_1),\dots,\Sigma(s_d) \in L_I$ are such $\que$-linearly independent vectors corresponding to these successive minima, \eqref{I_sm-1} follows immediately by combining~\eqref{mink_jarnik} with~\eqref{det_norm} and Lemma~\ref{ht_ineq}. To obtain~\eqref{I_sm-2}, notice that the second inequality of~\eqref{mink_jarnik} implies that
$$\mu(L_I) \leq \frac{d^{3/2}}{2} \lambda_d.$$
On the other hand, Lemma~\ref{ht_ineq} guarantees that $\lambda_1 \geq 1$, and so $\prod_{i=1}^d \lambda_i \geq \lambda_d$. Thus applying now the first inequality of~\eqref{mink_jarnik}, we see that
$$\mu(L_I) \leq \frac{d^{3/2}}{2} \lambda_d \leq \frac{d^{3/2}}{2} \prod_{i=1}^d \lambda_i \leq \frac{d^{3/2}}{2} \det L_I = \frac{d^{3/2}}{2} \Nn_K(I) \sqrt{ |\Delta_K| },$$
where the last equality is given by ~\eqref{det_norm}.
\endproof

\begin{lem} \label{I+_sm} There exist $\que$-linearly independent elements $\alpha_1,\dots,\alpha_d \in I^+$ such that
$$\prod_{i=1}^d h(\alpha_i)  \leq \left( 3d \sqrt{d} \right)^d \left( \Nn_K(I) \sqrt{ |\Delta_K| } \right)^{d+1}.$$
\end{lem}

\proof
A collection of elements $\alpha_1,\dots,\alpha_d$ is $\que$-linearly independent in $I^+$ if and only if the collection $\Sigma(\alpha_1),\dots,\Sigma(\alpha_d)$ is $\que$-linearly independent in $L_I^+$. Take a collection $s_1,\dots,s_d \in I$ as guaranteed by Lemma~\ref{I_sm} above. Notice that Lemma~\ref{ht_ineq} implies that $\mu(L_I) \geq 1/2$. Then combining~\eqref{smL+} with Lemma~\ref{ht_ineq}, we see that there exist such $\que$-linearly independent in $I^+$ so that
\begin{eqnarray*}
\prod_{i=1}^d h(\alpha_i) & \leq & \prod_{i=1}^d |\Sigma(\alpha_i)| \leq (2\mu(L_I)+1) \prod_{i=2}^d 2 |\Sigma(s_i)| (\mu(L_I) + 1) \\
& \leq & 6^d \mu(L_I)^d \prod_{i=2}^d |\Sigma(s_i)| \leq 6^d \mu(L_I)^d \prod_{i=1}^d |\Sigma(s_i)| \\
& \leq & \left( 3 d^{3/2} \right)^d \left( \Nn_K(I) \sqrt{ |\Delta_K| } \right)^{d+1},
\end{eqnarray*}
where the last inequality follows by Lemma~\ref{I_sm}.
\endproof

\begin{lem} \label{Ib_sm} Let $\bbb = \{ \beta_1,\dots,\beta_d \}$ be a positive basis for the ideal $I$ and $G(\bbb) = I^+ \setminus S(\bbb)$ be the corresponding set of gaps. For each $1 \leq i \leq d$, let $\beta'_i = \sum_{j=1, j \neq i}^d \beta_j$. Then there exist $\que$-linearly independent gaps $\alpha_1,\dots,\alpha_d \in G(\bbb)$ such that
$$h(\alpha_i) \leq \left( h(\beta_i/\beta'_i)^d + 1\right) h(\beta'_i)^d.$$
\end{lem}

\proof
Notice that $G(\Sigma(\bbb)) = \Sigma(G(\bbb))$. Then Lemma~\ref{gap_ht_1} implies existence of $\que$-linearly independent $\alpha_1,\dots,\alpha_d \in G(\bbb)$ such that 
\begin{eqnarray*}
|\Sigma(\alpha_i)| & = & \max_{1 \leq m \leq d}  \left\{ \left( \max_{1 \leq k \leq d} \left[ \frac{\sigma_k(\beta_i)}{\sum_{j=1, j \neq i}^d \sigma_k(\beta_j)} \right] + 1 \right) \sum_{j=1, j \neq i}^d \sigma_m(\beta_j) - \sigma_m(\beta_i)\right\} \\
& \leq & \max_{1 \leq m \leq d}  \left\{ \left( \max_{1 \leq k \leq d} \sigma_k \left( \frac{\beta_i}{\beta'_i} \right) + 1 \right) \sigma_m \left( \beta'_i  \right) - \sigma_m (\beta_i) \right\} \\
& \leq & \left( |\Sigma(\beta_i/\beta'_i)| + 1\right) |\Sigma(\beta'_i)| \leq \left( h(\beta_i/\beta'_i)^d + 1\right) h(\beta'_i)^d,
\end{eqnarray*}
by Lemma~\ref{ht_ineq}, and $h(\alpha_i) \leq |\Sigma(\alpha_i)|$ also by Lemma~\ref{ht_ineq}. This completes the proof.
\endproof

Now Theorem~\ref{ideal_gaps} follows by combining Lemmas~\ref{I_sm}, \ref{I+_sm} and \ref{Ib_sm}.

\bigskip

\noindent
{\bf Acknowledgement:} We would like to sincerely thank the anonymous referee whose excellent suggestions helped to significantly improve our paper.

\bibliographystyle{plain}  
\bibliography{semigroups}    
\end{document}